\documentclass[12pt]{amsart}
\usepackage[russian]{babel}
\usepackage{amsmath,amsfonts,amssymb,euscript}

\theoremstyle{plain}
\newtheorem{Theorem}{╥хюЁхьр}
\newtheorem{THEOREM}{╥хюЁхьр}

\newtheorem*{Corollary}{╤ыхфёЄтшх}

\newcommand{\intl}{\mathop{\int}\limits}

\newcommand{\suml}{\mathop{\sum}\limits}

\renewcommand\Re{\mathrm{Re}\,}
\renewcommand\Im{\mathrm{Im\,}}
\renewcommand{\le}{\leqslant}
\renewcommand{\ge}{\geqslant}

\newcommand\wt{\widetilde}

\def\al{\alpha}
\def\la{\lambda}

\def\eps{\varepsilon}

\def\bR{{\mathbb R}}
\def\bN{{\mathbb N}}
\def\bZ{{\mathbb Z}}

\newlength{\lenun}
\newlength{\lendu}
\def\s{{\mathbf s}}

\begin{document}

{\bf ╙─╩ 517.984.54}

\begin{center}{\bf\large  ┬юёёЄрэютыхэшх яюЄхэЎшрыр юяхЁрЄюЁр ╪ЄєЁьр--╦шєтшыы  яю ъюэхўэюьє эрсюЁє ёюсёЄтхээ√ї чэрўхэшщ ш эюЁьшЁютюўэ√ї ўшёхы} \end{center}
\medskip
\centerline{└. ╠. ╤ртўєъ\footnote{╨рсюЄр яюффхЁцрэр уЁрэЄюь ╨═╘ 14-01-00754.}}

\vskip 1cm

╠√ шчєўрхь ъырёёшўхёъє■ юсЁрЄэє■ чрфрўє тюёёЄрэютыхэш  тх∙хёЄтхээюую яюЄхэЎшрыр $q$ юяхЁрЄюЁр ╪ЄєЁьр--╦шєтшыы  $L$,
чрфрээюую т яЁюёЄЁрэёЄтх $L_2[0,\pi]$ яю ёяхъЄЁры№эющ ЇєэъЎшш $\rho_L(\la)$. ╠√ ЁрёёьрЄЁштрхь ъЁрхт√х єёыютш  ─шЁшїых,
їюЄ  эр° ьхЄюф ш тёх Ёхчєы№ЄрЄ√ ЁрсюЄ√ ыхуъю яхЁхэюё Єё  эр ёыєўрщ ъЁрхт√ї єёыютшщ ─шЁшїых--═хщьрэр, Є.х. єёыютшщ тшфр
$y(0)=y^{[1]}(\pi)=0$ (чфхё№ $y^{[1]}(x)=y'(x)-\sigma(x)y$, $\sigma=\int q(\xi)\,d\xi$). ╤ыєўрщ ъЁрхт√ї єёыютшщ ═хщьрэр
шьххЄ ётю■ ёяхЎшЇшъє (ёь. \cite{Amb} шыш \cite{Yu1}[├ыртр I, \S 2]) ш ЄЁхсєхЄ юЄфхы№эюую шчєўхэш . ═р° яюфїюф яючтюы хЄ
яюыєўрЄ№ Ёхчєы№ЄрЄ√ ш ъръ фы  ъырёёшўхёъшї яюЄхэЎшрыют, Єръ ш фы  яюЄхэЎшрыют--ЁрёяЁхфхыхэшщ --- ь√ яЁхфяюырурхь, ўЄю
$q$ яЁшэрфыхцшЄ яЁюёЄЁрэёЄтє ╤юсюыхттр  $W_2^\theta[0,\pi]$, $\theta\in(-1,\infty)$. ╦хуъю тшфхЄ№, ўЄю яюёъюы№ъє ь√
шчєўрхь юяхЁрЄюЁ эр ъюэхўэюь юЄЁхчъх, ёяхъЄЁры№эр  ЇєэъЎш  хёЄ№ ЇєэъЎш  ёърўъют. ╧Ёш ¤Єюь Єюўъш ЁрчЁ√тют ёютярфр■Є ё
ёюсёЄтхээ√ьш чэрўхэш ьш $\{\la_n\}_1^\infty$ юяхЁрЄюЁр $L_D$, р тхышўшэ√ ёърўъют Ёртэ√ $\al_n^{-1}$, уфх
$\{\al_n\}_1^\infty$ --- эрсюЁ эюЁьшЁютюўэ√ї ўшёхы юяхЁрЄюЁр $L_D$. ╥ръшь юсЁрчюь, т ърўхёЄтх ёяхъЄЁры№эющ шэЇюЁьрЎшш
ь√ ЁрёёьрЄЁштрхь фтх яюёыхфютрЄхы№эюёЄш $\{\la_n\}_1^\infty\cup\{\al_n\}_1^\infty$. ╒юЁю°ю шчтхёЄэю, ўЄю яю фтєь фрээ√ь
тх∙хёЄтхээ√ь яюёыхфютрЄхы№эюёЄ ь (яЁш єёыютшш яюыюцшЄхы№эюёЄш ўшёхы $\al_n$ ш яЁш т√яюыэхэшш юяЁхфхыхээ√ї
рёшьяЄюЄшўхёъшї ёююЄэю°хэшщ) яюЄхэЎшры $q$ тюёёЄрэртыштрхЄё  юфэючэрўэю. ╠√ шчєўрхь тюяЁюё ю тюёёЄрэютыхэшш яюЄхэЎшрыр
яю ъюэхўэюьє эрсюЁє ёяхъЄЁры№э√ї фрээ√ї. ┬ ърўхёЄтх Єръюую эрсюЁр ь√ шёяюы№чєхь ўшёыр $\{\la_n\}_1^N\cup\{\al_n\}_1^N$
ш эхёъюы№ъю фюяюыэшЄхы№э√ї ъю¤ЇЇшЎшхэЄют т рёшьяЄюЄшўхёъшї яЁхфёЄртыхэш ї фы  $\la_n$ ш $\al_n$ (эряЁшьхЁ, фы  ёыєўр 
$q\in L_2[0,\pi]$ ь√ яЁхфяюырурхь шчтхёЄэ√ь ёЁхфэхх чэрўхэшх яюЄхэЎшрыр). ▀ёэю, ўЄю Єрър  шэЇюЁьрЎш  эх яючтюы хЄ эрщЄш
яюЄхэЎшры $q$ т ЄюўэюёЄш, їюЄ  їюЁю°ю шчтхёЄхэ  тэ√щ рыуюЁшЄь, яючтюы ■∙шщ яюёЄЁюшЄ№ ЇєэъЎш■ $q_N$, яЁшсышцр■∙є■
яюЄхэЎшры. ╓хы№ ¤Єющ ЁрсюЄ√ --- фрЄ№ юЎхэъш эр эюЁьє ЁрчэюёЄш $\|q-q_N\|$ т яЁюёЄЁрэёЄтх ╤юсюыхтр $W_2^\tau[0,\pi]$ фы 
$\tau\in[-1,\theta]$. ╬Ўхэъш Єръюую тшфр с√ыш шчтхёЄэ√ фртэю, эю с√ыш ыюъры№э√ьш --- ёюфхЁцрыш яюёЄю ээ√х, чртшё ∙шх юЄ
яюЄхэЎшрыр $q$, ўЄю фхырыю шї ьрыю яЁшьхэшь√ьш фы  яЁръЄшўхёъюую шёяюы№чютрэш . ╤ фЁєующ ёЄюЁюэ√,  ёэю, ўЄю яЁш
юЄёєЄёЄтшш тююс∙х ъръющ--ышсю ряЁшюЁэющ шэЇюЁьрЎшш эюЁьр $\|q-q_N\|$ ьюцхЄ юърчрЄ№ё  ёъюы№ єуюфэю сюы№°ющ. ╠√
ЁрёёьюЄЁшь фтр тшфр ряЁшюЁэющ шэЇюЁьрЎшш: шэЇюЁьрЎш■ ю яюЄхэЎшрых шыш ю яюыэюь эрсюЁх ёяхъЄЁры№э√ї фрээ√ї.

▌Єр ёЄрЄ№  яЁюфюыцрхЄ шёёыхфютрэш  ртЄюЁют, эрўрЄ√х т ЁрсюЄрї \cite{SS1}--\cite{SS7}. ╠√ яЁштхфхь чфхё№ ыш°№ ъЁрЄъшщ
шёЄюЁшўхёъшщ юсчюЁ --- яюфЁюсэє■ шэЇюЁьрЎш■ ўшЄрЄхы№ ьюцхЄ эрщЄш т ьюэюуЁрЇш ї ╦хтшЄрэр \cite{Lev1}, ╠рЁўхэъю
\cite{Ma1}, ╫рфрэр ш ╤рсрЄ№х \cite{CS}, ╧ю°хы  ш ╥ЁєсютшЎр \cite{PT}, ╘Ёрщышэур ш ▐Ёъю \cite{FYu} ш \cite{Yu1} ш т
ёЄрЄ№ ї, ЎшЄшЁютрээ√ї эшцх. ╚ёёыхфютрэш  т юсырёЄш юсЁрЄэющ чрфрўш ╪ЄєЁьр--╦шєтшыы  тюёїюф Є ъ ЁрсюЄрь └ьсрЁЎєь эр
\cite{Amb} ш ┴юЁур \cite{Bo}. ┬ яюёыхфэхщ ЁрсюЄх с√ыр шёёыхфютрэр чрфрўр тюёёЄрэютыхэш  яюЄхэЎшрыр яю фтєь ёяхъЄЁрь.
╟рфрўр тюёёЄрэютыхэш  яюЄхэЎшрыр яю ёяхъЄЁры№эющ ЇєэъЎшш тючэшъыр т ЁрсюЄрї ├хы№Їрэфр ш ╦хтшЄрэр. ╬ЄьхЄшь, ўЄю ё Єюўъш
чЁхэш  тюяЁюёют ЁрчЁх°шьюёЄш, ¤Єш чрфрўш ¤ътштрыхэЄэ√ (ёь. \cite{Levi}) ш  ты ■Єё  ўрёЄэ√ьш ёыєўр ьш тюёёЄрэютыхэш 
яюЄхэЎшрыр яю ЇєэъЎшш ┬хщы --╥шЄўьрЁ°р (ёь. \cite{Yu1}[├ыртр I, \S 3-5]). ═х яЁхЄхэфє  эр яюыэюЄє, яЁштхфхь чфхё№
ёё√ыъш эр ЁрсюЄ√ ├хы№Їрэфр, ╦хтшЄрэр ш ├рё√ьютр \cite{GaL}, \cite{GeL}, \cite{Lev1}, ╠рЁўхэъю ш ╬ёЄЁютёъюую \cite{Ma},
\cite{Ma1} \cite{MO} ш ╩Ёхщэр \cite{Kr1}, \cite{Kr2}. ┬ ўрёЄэюёЄш, т ¤Єшї ЁрсюЄрї с√ыш эрщфхэ√ эхюсїюфшь√х ш
фюёЄрЄюўэ√х єёыютш  фы  ёє∙хёЄтютрэшш ш хфшэёЄтхээюёЄш Ёх°хэш  юсЁрЄэющ чрфрўш тюёёЄрэютыхэш  яюЄхэЎшрыр яю
ёяхъЄЁры№эющ ЇєэъЎшш фы  ёыєўр  $q\in W_2^n[0,\pi]$, $n\in\{0\}\cup\bN$. ╧Ё ь√х ш юсЁрЄэ√х ёяхъЄЁры№э√х чрфрўш фы 
юяхЁрЄюЁют ╪ЄєЁьр--╦шєтшыы  ё яюЄхэЎшрырьш --- ЁрёяЁхфхыхэш ьш $q\in W_2^{-1}[0,\pi]$ шчєўрышё№ т ЁрсюЄрї ртЄюЁют
\cite{SS1}--\cite{SS7}, ├Ёшэштр ш ╠шъшЄ■ър \cite{HM1}--\cite{HM5}, ─цръютр ш ╠шЄ ушэр \cite{DM1}, \cite{DM2}. ┬юяЁюё√
єёЄющўштюёЄш Ёх°хэш  юсЁрЄэющ чрфрўш (т ыюъры№эющ ЇюЁьх) шёёыхфютрышё№ т ЁрсюЄрї ╠рЁўхэъю ш ╠рёыютр \cite{MM}, ▐Ёъю
\cite{Yu}, ╨ сє°ъю \cite{R1}, \cite{R2}, ├Ёшэштр \cite{Hr} ш ╠ръ╦рЇышэ \cite{McL}. ╬ЄьхЄшь Єръцх сышчъшх яю Єхьх ЁрсюЄ√
╒шЄЁшър \cite{Hit} ш ╩юЁюЄ хтр \cite{Ko}. ╨ртэюьхЁэ√х юЎхэъш єёЄющўштюёЄш с√ыш яюыєўхэ√ ртЄюЁрьш т ЁрсюЄрї \cite{SS5} ш
\cite{SS6}. ┬ ўрёЄэюёЄш, Ёхчєы№ЄрЄ√ ЁрсюЄ√ \cite{SS6}  ты ■Єё  ъы■ўхт√ьш фы  ¤Єющ ёЄрЄ№ш. ┬юяЁюё√ тюёёЄрэютыхэш 
яюЄхэЎшрыр яю ъюэхўэюьє эрсюЁє ёяхъЄЁры№э√ї фрээ√ї ш тюяЁюё√ ўшёыхээюую ьюфхышЁютрэш  ¤Єющ чрфрўш шчєўрышё№ т ЁрсюЄрї
\cite{PT}, \cite{An}, \cite{BSKM}, \cite{Hald}, \cite{Paine}, \cite{Ro} ш \cite{RS}. ╬ЄьхЄшь Єръцх ЁрсюЄє ╠рЁыхЄЄр ш
┬рщърЁфр \cite{MW}, т ъюЄюЁющ с√ыш фюърчрэ√ юЎхэъш эр эюЁьє $\|q-q_N\|$, сышчъшх ъ фюърчрээ√ь чфхё№. ┬ эрёЄю ∙хщ ёЄрЄ№х
(р Єръцх т ЁрсюЄх \cite{SS7}, уфх с√ы ЁрчюсЁрэ ёыєўрщ юсЁрЄэющ чрфрўш яю фтєь ёяхъЄЁрь) ь√ яюыєўрхь  тэ√х юЎхэъш эр
ёъюЁюёЄ№ ёїюфшьюёЄш $q_N\to q$ тю тёхщ °ърых яЁюёЄЁрэёЄт ╤юсюыхтр ё ЁртэюьхЁэ√ьш ъюэёЄрэЄрьш.

─ы  єфюсёЄтр ўшЄрЄхы  ь√ эрўэхь ё ЇюЁьєышЁютъш Ёхчєы№ЄрЄют фы  ъырёёшўхёъюую ёыєўр  $q\in L_2[0,\pi]$ --- т ¤Єюь ёыєўрх
эрь эх яюЄЁхсєхЄё  эшъръшї фюяюыэшЄхы№э√ї яюёЄЁюхэшщ. ┬ю тЄюЁюь ярЁруЁрЇх ЁрсюЄ√ ь√ фюърцхь ¤Єш Ёхчєы№ЄрЄ√ фы  юс∙хую
ёыєўр  $q\in W_2^\theta[0,\pi]$, $\theta\ge-1$.

\section{╤ыєўрщ $q\in L_2[0,\pi]$.}

\vskip 8pt


╨рёёьюЄЁшь юяхЁрЄюЁ ╪ЄєЁьр--╦шєтшыы , яюЁюцфхээ√щ т яЁюёЄЁрэёЄтх $L_2[0,\pi]$ фшЇЇхЁхэЎшры№э√ь т√Ёрцхэшхь
$$
l(y) = -y'' + q(x) y
$$
ё яюЄхэЎшрыюь  $q\in L_2[0,\pi]$. ┬ ¤Єющ ЁрсюЄх ь√ юуЁрэшўшьё  ёыєўрхь ъЁрхт√ї єёыютшщ ─шЁшїых\footnote{╟рьхЄшь, ўЄю
яюыєўхээ√х чфхё№ Ёхчєы№ЄрЄ√ ыхуъю яхЁхэюё Єё  ш эр ёыєўрщ ъЁрхт√ї єёыютшщ $y(0) =y'(\pi)-hy(\pi)=0$.}, юсючэрўшт ўхЁхч
$L_D$ юяхЁрЄюЁ, фхщёЄтє■∙шщ яю яЁртшыє $L_Dy=l(y)$ эр юсырёЄш юяЁхфхыхэш 
$$
\mathcal D(L_D) =\{y \in W^2_2[0,1]\vert \ \, y(0) =y(\pi) =0\}.
$$
╫хЁхч  $W^\theta_2[0,\pi]$ чфхё№ ш фрыхх ь√ юсючэрўрхь яЁюёЄЁрэёЄтр ╤юсюыхтр (ярЁрьхЄЁ $\theta\in[-1,+\infty)$). ╒юЁю°ю
шчтхёЄэю, ўЄю юяхЁрЄюЁ $L_D$ шьххЄ фшёъЁхЄэ√щ ёяхъЄЁ. ╠√ юсючэрўшь ўхЁхч $\{\la_k\}_1^\infty$ хую ёюсёЄтхээ√х чэрўхэш ,
яЁшўхь эєьхЁрЎш■ сєфхь тхёЄш т яюЁ фъх тючЁрёЄрэшщ ьюфєы  ш ё єўхЄюь рыухсЁршўхёъющ ъЁрЄэюёЄш. ╧Ёш Єръющ эєьхЁрЎшш
ёяЁртхфышт√ рёшьяЄюЄшўхёъшх яЁхфёЄртыхэш  (ёь. \cite{Lev1} шыш \cite{SS4})
\begin{equation} \label{eq:spas}
\sqrt{\lambda_k}=k+\dfrac{q_0}{2k}+ \dfrac{b_k}{k},
\end{equation}
уфх  юёЄрЄъш $\{b_k\}_1^{\infty}\in\ell_2$, р
\begin{equation}\label{eq:spas1}
q_0 = \frac{\sigma(\pi)}{\pi}=\frac1{\pi}\int_0^\pi q(x)\,dx.
\end{equation}
╟фхё№  $\sigma$ --- яхЁтююсЁрчэр  яюЄхэЎшрыр $q$, єфютыхЄтюЁ ■∙р  єёыютш■ $\sigma(0)=0$ (єфюсёЄтю т√сюЁр шьхээю Єръюую
єёыютш  сєфхЄ юс· ёэхэю тю тЄюЁюь ярЁруЁрЇх). ╟рьхЄшь, ўЄю ¤Єш рёшьяЄюЄшўхёъшх ЇюЁьєы√  ¤ътштрыхэЄэ√
\begin{equation*}
\lambda_k =  k^2 + q_0 + \wt{b}_{k},\qquad\text{уфх}\ \ \sum_1^\infty |\wt b_k|^2 <\infty.
\end{equation*}

┬ ёыєўрх тх∙хёЄтхээюую яюЄхэЎшрыр тёх ёюсёЄтхээ√х ўшёыр $\{\la_k\}_1^\infty$ тх∙хёЄтхээ√ ш яЁюёЄ√. ┬ юс∙хь ёыєўрх ¤Єш
ўшёыр ыхцрЄ тэєЄЁш эхъюЄюЁющ ярЁрсюы√ $\Re\la=k\cdot\Im\la^2+\la_0$, уфх $k>0$, р $\la_0\in\bR$. ─юуютюЁшьё , ўЄю т
\eqref{eq:spas} ш тхчфх фрыхх тхЄт№ ъюЁэ  т√сшЁрхЄё  єёыютшхь $\rm{Arg}\,\sqrt{\la}\in(-\pi,\pi]$. ╬сючэрўшь ЄхяхЁ№
ўхЁхч $s(x,\la)$ Ёх°хэшх єЁртэхэш  $l(y)=\la y$, єфютыхЄтюЁ ■∙хх эрўры№э√ь єёыютш ь $s(0,\la)=0$,
$s'(0,\la)=\sqrt{\la}$. ╦хуъю тшфхЄ№, ўЄю эєыш ЇєэъЎшш $s(\pi,\la)$ ёютярфр■Є ё ўшёырьш $\{\la_k\}_1^\infty$. ┬ ёыєўрх
тх∙хёЄтхээюую яюЄхэЎшрыр ўшёыр
\begin{equation}\label{eq:nc}
\al_k:=\begin{cases}\displaystyle{\int_0^\pi s^2(x,\la_k)\,dx,}\qquad &\text{хёыш}\ \la_k\ne0;\\
\displaystyle{\lim_{\la\to0}\frac1{\la}\int_0^\pi s^2(x,\la)\,dx,}\qquad &\text{хёыш}\ \la_k=0,\end{cases}\qquad k\ge1,
\end{equation}
эрч√тр■Є \textit{эюЁьшЁютюўэ√ьш ъюэёЄрэЄрьш}. ╠√ ёюїЁрэшь ¤Єю эрчтрэшх ш фы  ёыєўр , ъюуфр яюЄхэЎшры яЁшэшьрхЄ
ъюьяыхъёэ√х чэрўхэш . ╚чтхёЄэю (ёь. \cite{Lev1} шыш \cite{SS4}), ўЄю
\begin{equation}\label{eq:ncas}
\al_k=\frac{\pi}2+\dfrac{a_k}k,\qquad\text{уфх}\ \ \{a_k\}_1^\infty\in\ell_2.
\end{equation}
╟рьхЄшь, ўЄю т ъырёёшўхёъшї ЁрсюЄрї эюЁьшЁютюўэ√х ўшёыр ёЄЁю Єё  яю ЇєэъЎшш $s(x,\la)$ ё эрўры№э√ьш єёыютш ьш
$s(0,\la)=0$, $s'(0,\la)=1$. ╦хуъю тшфхЄ№, ўЄю юфшэ ёяюёюс ётюфшЄё  ъ фЁєуюьє єьэюцхэшхь эр яюёыхфютрЄхы№эюёЄ№
$\{\sqrt{\la_k}\}_1^\infty$. ╬с·хфшэхэшх фтєї яюёыхфютрЄхы№эюёЄхщ
$$
\{\la_k\}_1^\infty\cup\{\al_k\}_1^\infty
$$
ь√ эрчютхь \textit{ёяхъЄЁры№э√ьш фрээ√ьш}, юЄтхўр■∙шьш яюЄхэЎшрыє $q$. ╬ЄьхЄшь, ўЄю т ёыєўрх тх∙хёЄтхээюую яюЄхэЎшрыр
чрфрэшх эрсюЁр ёяхъЄЁры№э√ї фрээ√ї ¤ътштрыхэЄэю чрфрэш■ ёяхъЄЁры№эющ ЇєэъЎшш юяхЁрЄюЁр $L_D$, юяЁхфхы хьющ ёєььющ
$$
\rho(\la)=\suml_{\la_k<\la}\frac{1}{\al_k}.
$$

╒юЁю°ю шчтхёЄэю (ёь., эряЁшьхЁ, \cite{Ma}), ўЄю ёяхъЄЁры№э√х фрээ√х юфэючэрўэю юяЁхфхы ■Є яюЄхэЎшры, т ўрёЄэюёЄш,
$q_0=\lim_{k\to\infty}(\la_k-k^2)$. 
╠√ ЁрёёьюЄЁшь тюяЁюё ю тюёёЄрэютыхэшш яюЄхэЎшрыр яю \textit{ъюэхўэюьє эрсюЁє ёяхъЄЁры№э√ї фрээ√ї}
\begin{equation}\label{eq:findata}
\{q_0\}\cup\{\la_k\}_1^N\cup\{\al_k\}_1^N,
\end{equation}
уфх ўшёыю $q_0\in\bR$ яЁюшчтюы№эю,  тх∙хёЄтхээ√х ўшёыр $\la_k$ єфютыхЄтюЁ ■Є ёююЄэю°хэш ь
$$
\la_1<\la_2<\dots<\la_N,
$$
р тёх $\al_k$ яюыюцшЄхы№э√. ╟рьхЄшь, ўЄю хёыш $q=q_0=\text{const}$, Єю
$$
\la_k=\la_k^0=k^2+q_0,\qquad\text{р}\ \ \al_k=\al_k^0=\frac{\pi}{2}+\frac{\pi q_0}{2k^2}.
$$
╧єёЄ№ ЄхяхЁ№ $q\in L_2[0,\pi]$ --- эхъюЄюЁ√щ тх∙хёЄтхээючэрўэ√щ яюЄхэЎшры, $\{\la_k\}_1^\infty\cup\{\al_k\}_1^\infty$
--- хую ёяхъЄЁры№э√х фрээ√х, яЁшўхь эрь шчтхёЄхэ Єюы№ъю ъюэхўэ√щ эрсюЁ ёяхъЄЁры№э√ї фрээ√ї \eqref{eq:findata}. ╤юуырёэю
ъырёёшўхёъющ ЄхюЁшш юсЁрЄэющ чрфрўш (ёь., эряЁшьхЁ, \cite{MO} шыш \cite{Ma}), ёє∙хёЄтєхЄ хфшэёЄтхээ√щ яюЄхэЎшры $q_N$,
ёяхъЄЁры№э√х фрээ√х ъюЄюЁюую шьх■Є тшф
$$
\{\la_1,\dots,\la_N,\la_{N+1}^0,\la_{N+2}^0,\dots\}\cup\{\al_1,\dots,\al_N,\al_{N+1}^0,\al_{N+2}^0,\dots\},
$$
уфх $\la_k^0$ ш $\al_k^0$ --- ёюсёЄтхээ√х ўшёыр ш, ёююЄтхЄёЄтхээю, эюЁьшЁютюўэ√х ўшёыр юяхЁрЄюЁр ё яюЄхэЎшрыюь
$q_0=\frac1{\pi}\int_0^\pi q(x)\,dx$. ┴юыхх Єюую, ¤ЄюЄ яюЄхэЎшры ьюцхЄ с√Є№ эрщфхэ т ърўхёЄтх Ёх°хэш  ёшёЄхь√ $2N$
ышэхщэ√ї єЁртэхэшщ (Єръющ тшф яЁшэшьрхЄ чфхё№ єЁртэхэшх ├хы№Їрэфр--╦хтшЄрэр, ёь., эряЁшьхЁ \cite{Yu}). ▌Єє ёшёЄхьє
ьюцэю ышсю Ёх°рЄ№ эхяюёЁхфёЄтхээю, ышсю яЁшьхэ   шчтхёЄэє■ $2N$--°руютє■ яЁюЎхфєЁє (ёь. \cite{SS6} шыш \cite{PT}).
╬яЁхфхыхээ√щ Єръшь юсЁрчюь яюЄхэЎшры $q_N$ ь√ сєфхь эрч√трЄ№ \textit{$2N$--ряяЁюъёшьрЎшхщ яюЄхэЎшрыр $q$, яюёЄЁюхээющ
яю ъюэхўэюьє эрсюЁє ёяхъЄЁры№э√ї фрээ√ї \eqref{eq:findata}}. ┼ёЄхёЄтхээю юцшфрЄ№, ўЄю $q_N\to q$  яю эюЁьх яЁюёЄЁрэёЄтр
$L_2[0,\pi]$ яЁш $N\to\infty$. ▌Єю фхщёЄтшЄхы№эю Єръ, р шьхээю, ёяЁртхфыштю ёыхфє■∙хх єЄтхЁцфхэшх.

\begin{Theorem}\label{tm:1} ╧єёЄ№ $q\in L_2[0,\pi]$ --- яЁюшчтюы№эр  тх∙хёЄтхээючэрўэр  ЇєэъЎш , $\{\la_k\}_1^\infty$
--- ёюсёЄтхээ√х чэрўхэш , р $\{\al_k\}_1^\infty$ --- эюЁьшЁютюўэ√х ўшёыр юяхЁрЄюЁр $L_D$ ё яюЄхэЎшрыюь $q$.
╧Ёхфяюыюцшь, ўЄю эрь шчтхёЄхэ ъюэхўэ√щ эрсюЁ ёяхъЄЁры№э√ї фрээ√ї \eqref{eq:findata}. ╧єёЄ№  $q_N$ ---
$2N$-ряяЁюъёшьрЎш  яюЄхэЎшрыр $q$, яюёЄЁюхээр  яю ¤Єюьє эрсюЁє. ╥юуфр
$$
\|q_N - q\|_{L_2} \leqslant C(q)\left(\sum_{k=N+1}^\infty |a_k|^2+|b_k|^2\right)^{1/2} \to 0, \quad \text{яЁш}\ \,
N\to\infty,
$$
уфх $a_k$ ш $b_k$ юяЁхфхы ■Єё  яюЄхэЎшрыюь $q$ т \eqref{eq:spas1} ш \eqref{eq:ncas} ёююЄтхЄёЄтхээю, р ўшёыю $C=C(q)$
Єръцх чртшёшЄ юЄ $q$. ┴юыхх Єюую, ъюэёЄрэЄр $C$ ьюцхЄ с√Є№ т√сЁрэр хфшэ√ь юсЁрчюь т °рЁх $\|q\|_{L_2}\le R$, Є.х.
\begin{equation}\label{eq:appr0}
\|q_N - q\|_{L_2} \leqslant C(R)\left(\sum_{k=N+1}^\infty |a_k|^2+|b_k|^2\right)^{1/2} \to 0, \quad \text{яЁш}\ \,
N\to\infty,
\end{equation}
уфх $C$ чртшёшЄ Єюы№ъю юЄ $R$.
\end{Theorem}
╧хЁтюх єЄтхЁцфхэшх ¤Єющ ЄхюЁхь√ (фы  ёыєўр  $q_0=0$) ёыхфєхЄ шч Ёхчєы№ЄрЄют ЁрсюЄ \cite{MO} ш \cite{R1}, \cite{R2}.
┬ЄюЁюх єЄтхЁцфхэшх ёыхфєхЄ шч ЁрсюЄ√ ртЄюЁют \cite{SS6} (ь√ фюърцхь хую тю тЄюЁюь ярЁруЁрЇх т сюыхх юс∙хщ ёшЄєрЎшш).

╤ Єюўъш чЁхэш  яЁръЄшўхёъюую шёяюы№чютрэш  ЄхюЁхьр \ref{tm:1} схёяюыхчэр т юЄёєЄёЄтшш ряЁшюЁэющ шэЇюЁьрЎшш, яючтюы ■∙хщ
юЎхэшЄ№ ёъюЁюёЄ№ ёїюфшьюёЄш $q_N\to q$. ╟ртшёшьюёЄ№ ъюэёЄрэЄ√ $C$ юЄ эюЁь√ $\|q\|_{L_2}$ яючтюы хЄ эрфх Єё , ўЄю т
ърўхёЄтх ряЁшюЁэющ шэЇюЁьрЎшш ьюцэю шёяюы№чютрЄ№ єёыютшх $\|q\|_{L_2}\le R$. ╬ърч√трхЄё , юфэръю, ўЄю Єръюх єёыютшх эх
яючтюы хЄ юЎхэшЄ№ ёъюЁюёЄ№ ёїюфшьюёЄш Ё фр $\sum_1^\infty|a_k|^2+|b_k|^2$. ╫шёыр $a_k$ ш $b_k$  ты ■Єё  эхышэхщэ√ьш
ЇєэъЎшюэрырьш, чртшё ∙шьш юЄ яюЄхэЎшрыр, яЁшўхь шї ышэхщэ√х ўрёЄш --- ¤Єю ъю¤ЇЇшЎшхэЄ√ ╘єЁ№х ЇєэъЎшщ $q(x)$ ш
$(\pi-x)q(x)$ ёююЄтхЄёЄтхээю. ╥ръшь юсЁрчюь, ышсю ряЁшюЁэр  шэЇюЁьрЎш  фюыцэр ёюфхЁцрЄ№ фтр єёыютш : $\|q\|_{L_2}\le R$
ш $\sum_{k=N+1}^\infty|a_k|^2+|b_k|^2<\eps$, ышсю эрь эєцэр шэЇюЁьрЎш  ю ёъюЁюёЄш єс√трэш  ъю¤ЇЇшЎшхэЄют ╘єЁ№х
яюЄхэЎшрыр. ╥юуфр яюфїюф ∙хщ ряЁшюЁэющ шэЇюЁьрЎшхщ  ты хЄё  єёыютшх $\|q\|_\theta\le R$, $\theta>0$ (чфхё№ ўхЁхч
$\|\cdot\|_\theta$ ь√ юсючэрўрхь эюЁьє т яЁюёЄЁрэёЄтх ╤юсюыхтр $W_2^\theta[0,\pi]$). ┼∙х юфэр тючьюцэюёЄ№ юЎхэшЄ№
ёъюЁюёЄ№ ёЄЁхьыхэш  $q_N\to q$ --- юёырсшЄ№ эюЁьє т ыхтющ ўрёЄш \eqref{eq:appr0}, чрьхэшт эюЁьє $\|\cdot\|_{L_2}$ эр
$\|\cdot\|_{\tau-1}$, $\tau\in[0,1)$.

\begin{Theorem}\label{tm:2} ╧єёЄ№ $q\in W_2^\theta[0,\pi]$, $\theta\in[0,1/2)$, --- яЁюшчтюы№эр  тх∙хёЄтхээючэрўэр  ЇєэъЎш ,
яЁшўхь $\|q\|_\theta\le R$. ╧єёЄ№ $\{\la_k\}_1^\infty$
--- ёюсёЄтхээ√х чэрўхэш , р $\{\al_k\}_1^\infty$ --- эюЁьшЁютюўэ√х ўшёыр юяхЁрЄюЁр $L_D$ ё яюЄхэЎшрыюь $q$.
╧Ёхфяюыюцшь, ўЄю эрь шчтхёЄхэ ъюэхўэ√щ эрсюЁ ёяхъЄЁры№э√ї фрээ√ї \eqref{eq:findata}. ╧єёЄ№  $q_N$ ---
$2N$-ряяЁюъёшьрЎш  яюЄхэЎшрыр $q$, яюёЄЁюхээр  яю ¤Єюьє эрсюЁє. ╬сючэрўшь $\sigma(x):=\int_0^x q(t)\,dt$ ш
$\sigma_N(x):=\int_0^x q_N(t)\,dt$. ╥юуфр фы  ы■сюую $\tau\in[0,1]$, $\tau-1<\theta$, ёяЁртхфыштр юЎхэър
\begin{equation}\label{eq:rate1}
\|q_N -q\|_{\tau -1} \leqslant \|\sigma_N(x)-\sigma(x)\|_\tau \leqslant \frac{C(R)}{N^{1+\theta-\tau}},
\end{equation}
уфх $C(R)$ чртшёшЄ Єюы№ъю юЄ $R$, $\theta$ ш $\tau$. ┬ ўрёЄэюёЄш, яЁш  $\tau \in (1/2, 1]$ т√яюыэхэю
\begin{equation}\label{eq:rate2}
\left|\int_0^x (q_N(t) -q(t))\, dt\right| =|\sigma_N(x)-\sigma(x)| \leqslant \frac{C(R)}{N^{1+\theta-\tau}},
 \end{equation}
уфх $C(R)$ чртшёшЄ Єюы№ъю юЄ  $R$, $\theta$ ш $\tau$.
\end{Theorem}
╠√ фюърцхь ¤Єє ЄхюЁхьє т сюыхх юс∙хщ ёшЄєрЎшш тю тЄюЁюь ярЁруЁрЇх ЁрсюЄ√. ╟фхё№ ыш°№ юЄьхЄшь, ўЄю юЎхэър
\eqref{eq:rate2} ёыхфєхЄ эхяюёЁхфёЄтхээю шч \eqref{eq:rate1} т ёшыє эхяЁхЁ√тэюую тыюцхэш  $W^\tau_2 \hookrightarrow
C[0,1]$ ё юЎхэъющ эюЁь√ $ \|\cdot\|_C \leqslant C(\tau)\|\cdot\|_\tau$ фы  тёхї $\tau > 1/2$ (ёь, эряЁшьхЁ,
\cite[Ch.5]{BIN})).


╥хюЁхьр \ref{tm:2} юёЄрты хЄ тюяЁюё ю ёъюЁюёЄш ЁртэюёїюфшьюёЄш юЄъЁ√Є√ь, яюёъюы№ъє т яЁртющ ўрёЄш \eqref{eq:rate1}
ЇшуєЁшЁєхЄ эхшчтхёЄэр  тхыхўшэр $C(R)$. └ъъєЁрЄэ√щ рэрышч ЁрсюЄ√ \cite{Hr} яюърч√трхЄ, ўЄю $C(R)\le A\exp(aR)$ ё
эхъюЄюЁ√ьш рсёюы■Єэ√ьш ъюэёЄрэЄрьш $A$ ш $a$. ╬фэръю шёяюы№чє  ьхЄюф ЁрсюЄ√ \cite{SS2}, ¤Єє юЎхэъє ьюцэю ёє∙хёЄтхээю
єЄюўэшЄ№ --- ртЄюЁ√ яЁхфяюырур■Є яюёт ЄшЄ№ ¤Єюьє юЄфхы№эє■ ЁрсюЄє.

╧Ёхфяюыюцшь ЄхяхЁ№, ўЄю є эрё эхЄ эшъръющ ряЁшюЁэющ шэЇюЁьрЎшш ю яюЄхэЎшрых $q$, эю хёЄ№ шэЇюЁьрЎш  юс рёшьяЄюЄшўхёъюь
яютхфхэшш яюёыхфютрЄхы№эюёЄхщ $\{\al_k\}_1^\infty$ ш $\{\la_k\}_1^\infty$.

\begin{Theorem}\label{tm:3}
╧Ёхфяюыюцшь, ўЄю эрь шчтхёЄхэ ъюэхўэ√щ эрсюЁ ёяхъЄЁры№э√ї фрээ√ї \eqref{eq:findata}. ╧єёЄ№ фюяюыэшЄхы№эю шчтхёЄэю, ўЄю
фы  яюыэюую эрсюЁр ёяхъЄЁры№э√ї фрээ√ї $\{\la_n\}_1^\infty\cup\{\al_n\}_1^\infty$ т√яюыэхэ√ ёююЄэю°хэш  \eqref{eq:spas}
ш \eqref{eq:ncas}, яЁшўхь
\begin{equation*}
\sum_{k=1}^\infty(|a_k|^2+|b_k|^2) \leqslant r,  \quad \la_{k+1}\ge \la_k+ h,\quad\la_k\ge k^2,\quad\text{ш}\quad
\al_k\ge h,\ k\in\bN.
 \end{equation*}
╧єёЄ№  $q_N$ --- $2N$-ряяЁюъёшьрЎш  эхшчтхёЄэюую яюЄхэЎшрыр $q$, яюёЄЁюхээр  яю ¤Єюьє эрсюЁє. ╥юуфр
$$
\|q_N - q\|_{L_2} \leqslant C(r,h)\left(\sum_{k=N+1}^\infty |a_k|^2+|b_k|^2\right)^{1/2}.
$$
\end{Theorem}
╒рЁръЄхЁ чртшёшьюёЄш ъюэёЄрэЄ√ $C$ юЄ $r$ ш $d$ Єръцх ЄЁхсєхЄ Є∙рЄхы№эюую шчєўхэш  ш чфхё№ ь√ ¤ЄюЄ тюяЁюё юяєёърхь.

═р яЁръЄшъх яЁш Ёх°хэшш юсЁрЄэющ чрфрўш ёюсёЄтхээ√х чэрўхэш  $\la_k$ ш эюЁьшЁютюўэ√х ўшёыр $\al_k$, $k=1,\,\dots,N$
шчтхёЄэ√ эрь ыш°№ яЁшсышцхээю. ╧Ёхфяюыюцшь, ўЄю ърцфюх шч ¤Єшї ўшёхы шчьхЁхэю ё эхъюЄюЁющ ЄюўэюёЄ№■ $\varepsilon$, р
шьхээю,
\begin{equation}\label{eq:err1}
|\al_k -\wt \al_k|\leqslant\varepsilon, \qquad |\lambda_k^{1/2} -\wt \lambda_k^{1/2}| \leqslant\frac{\varepsilon}{k},
\quad k=1,2,\dots, N,
\end{equation}
уфх  $\{\widetilde\al_k\}_1^N$ ш   $\{\widetilde\lambda_k\}_1^N$  --- фрээ√х шчьхЁхэшщ, р $\{\al_k\}_1^N$ ш
$\{\la_k\}_1^N$ --- Єюўэ√х чэрўхэш . ─ы  єяЁю∙хэш  фры№эхщ°шї ЇюЁьєышЁютюъ ь√ сєфхь ёўшЄрЄ№, ўЄю $\eps$ фюёЄрЄюўэю
ьрыю, р шьхээю, ўЄю $\eps<1/e$. ─рыхх ь√ сєфхь юсючэрўрЄ№ ўхЁхч $q_N$ ряяЁюъёшьрЎш■ яюЄхэЎЎшрыр, яюёЄЁюхээє■ яю
ъюэхўэюьє эрсюЁє ёяхъЄЁры№э√ї фрээ√ї \eqref{eq:findata}, р ўхЁхч $\widetilde q_N$ юсючэрўшь $2N$--ряяЁюъёшьрЎш■,
яюёЄЁюхээє■ яю {\it ъюэхўэюьє эрсюЁє яЁшсышцхээ√ї ёяхъЄЁры№э√ї фрэээ√ї}
\begin{equation}\label{eq:findataerr}
\{q_0\},\{\widetilde\al_k\}_1^N\cup\{\widetilde\lambda_k\}_1^N,\qquad q_0\in\bR,\ \wt\al_k>0,\
\wt\la_1<\wt\la_2<\dots<\wt\la_N.
\end{equation}
▀ёэю, ўЄю $\wt q_N\nrightarrow q$ яЁш $N\to\infty$. ╥хь эх ьхэхх, яюуЁх°эюёЄ№ $\|\wt q_N-q\|$ ьюцэю юЎхэшЄ№.

\begin{Theorem}\label{tm:4}    ╧єёЄ№ $q\in W_2^\theta[0,\pi]$, $\theta\in[0,1/2)$, --- яЁюшчтюы№эр  тх∙хёЄтхээючэрўэр  ЇєэъЎш ,
яЁшўхь $\|q\|_\theta\le R$. ╧єёЄ№ $\{\la_k\}_1^\infty$
--- ёюсёЄтхээ√х чэрўхэш , р $\{\al_k\}_1^\infty$ --- эюЁьшЁютюўэ√х ўшёыр юяхЁрЄюЁр $L_D$ ё яюЄхэЎшрыюь $q$.
╧Ёхфяюыюцшь, ўЄю эрь шчтхёЄхэ ъюэхўэ√щ эрсюЁ яЁшсышцхээ√ї ёяхъЄЁры№э√ї фрээ√ї \eqref{eq:findataerr}. ╧єёЄ№  $\wt q_N$
--- $2N$-ряяЁюъёшьрЎш  яюЄхэЎшрыр $q$, яюёЄЁюхээр  яю ¤Єюьє эрсюЁє. ╥юуфр фы  яЁюшчтюы№эюую $\tau\in[0,1]$,
$\tau-1<\theta$, т√яюыэхэю
\begin{align*}
&\|\wt q_N -q\|_{\tau -1} \leqslant \|\wt \sigma_N(x)-\sigma(x)\|_\tau \leqslant
C(R)\left(\frac{\varepsilon}{1-2\tau} +\frac 1{N^{1+\theta-\tau}}\right)&\text{яЁш}\ \tau\in[0,1/2);\notag\\
&\|\wt q_N -q\|_{-1/2} \leqslant \|\wt \sigma_N(x)-\sigma(x)\|_{1/2} \leqslant
C(R)\left(\varepsilon(\ln N)^{1/2} +\frac 1{N^{1/2+\theta}}\right)&\text{яЁш}\ \tau=\frac12;\notag\\
&\|\wt q_N -q\|_{\tau-1} \leqslant \|\wt \sigma_N(x)-\sigma(x)\|_{\tau} \leqslant C(R)\left(\varepsilon N^{2\tau-1}
+\frac 1{N^{1+\theta-\tau}}\right)&\text{яЁш}\ \tau\in(1/2,1],\notag
\end{align*}
уфх $C(R)$ чртшёшЄ Єюы№ъю юЄ $R$, $\sigma(x):=\int_0^x q(t)\,dt$, р $\wt \sigma_N(x):=\int_0^x \wt q_N(t)\,dt$.
\end{Theorem}

╥хяхЁ№ яюёЄртшь хёЄхёЄтхээє■ чрфрўє: эрщЄш эюьхЁ $N$ фы  ъюЄюЁюую $2N$--ряяЁюъёшьрЎш  $\wt q_N$ эршыєў°шь юсЁрчюь
яЁшсышцрхЄ яюЄхэЎшры $q$ ш юЎхэшЄ№ тхышўшэє $\|q-\wt q_N\|$.
\begin{Corollary}  ┬ єёыютш ї ЄхюЁхь√ \ref{tm:4}
\begin{align*}
  & \|\wt \sigma_N(x)-\sigma(x)\|_\tau \leqslant
C(R)\frac{\varepsilon}{1-2\tau} &&\text{яЁш}\ N\ge\left(\frac{1-2\tau}{\varepsilon}\right)^{\frac 1{1+\theta-\tau}},\ &\tau\in[0,1/2);\\
  &\|\wt\sigma_N(x)-\sigma(x)\|_{1/2}\leqslant C(R)\varepsilon |\ln \varepsilon|^{1/2}
  &&\text{яЁш}\ N=\eps^{-\frac2{1+\theta}}, \ &\tau=\frac12;\\
&\|\wt\sigma_N(x)-\sigma(x)\|_{\tau}\leqslant C(R)\varepsilon^{\frac{1+\theta -\tau}{\theta+\tau}} &&\text{яЁш}\ \ N=
\varepsilon^{-\frac1{\theta+\tau}},\ &\tau\in(1/2,1].
\end{align*}
\end{Corollary}

\section{╬ёэютэющ Ёхчєы№ЄрЄ}

═ряюьэшь юяЁхфхыхэшх юяхЁрЄюЁр ╪ЄєЁьр--╦шєтшыы  ё яюЄхэЎшрыюь $q\in W_2^{\theta-1}[0,\pi]$, $\theta\ge0$. ╧єёЄ№
$\sigma(x)$
--- яЁюшчтюы№эр  юсюс∙хээр  яхЁтююсЁрчэр   ЇєэъЎшш $q$ (шчтхёЄэю, ёь., эряЁшьхЁ, \cite[├ы.1]{GeS}, ўЄю юэр юяЁхфхыхэр
ё ЄюўэюёЄ№■ фю ъюэёЄрэЄ√ ш ыхцшЄ т яЁюёЄЁрэёЄтх $W_2^\theta[0,\pi]$). ╧юыюцшь
$$
y^{[1]}(x) = y'(x) - \sigma(x) y(x).
$$
╥хяхЁ№ юяхЁрЄюЁ $L_D$ ь√ ьюцхь юяЁхфхышЄ№ т√Ёрцхэшхь
$$
L y = -\left( y^{[1]}\right)' -\sigma(x) y^{[1]} - \sigma^2(x) y
$$
эр юсырёЄш
$$
\mathcal D (L_D) =\{ y, y^{[1]} \in W^1_1[0,1]\ \vert    \ Ly\in L_2[0,1],\ \,  y(0)=y(\pi) =0 \}.
$$
▌ЄюЄ юяхЁрЄюЁ  ты хЄё  ЇЁхфуюы№ьют√ь юяхЁрЄюЁюь ё эєыхт√ьш шэфхъёрьш (хёыш $q$ тх∙хёЄтхээючэрўэр, Єю $L_D$
ёрьюёюяЁ цхэ) ш фшёъЁхЄэ√ь ёяхъЄЁюь. ╧юфЁюсэ√щ юсчюЁ ётющёЄт юяхЁрЄюЁр $L_D$ ўшЄрЄхы№ ьюцхЄ эрщЄш т ЁрсюЄрї \cite{SS1}
ш \cite{SS2}. ╤Ёрчє цх юЄьхЄшь, ўЄю фюсртыхэшх ъюэёЄрэЄ√ ъ ЇєэъЎшш $\sigma$ эх ьхэ хЄ юяхЁрЄюЁ $L_D$. ╧ю¤Єюьє фрыхх ь√
сєфхь ёўшЄрЄ№, ўЄю ЇєэъЎш  $\sigma$ яЁшэрфыхцшЄ ЇръЄюЁ--яЁюёЄЁрэёЄтє $W_2^\theta[0,\pi]/\{1\}$ (ёъры Ёэюх яЁюшчтхфхэшх
ш эюЁьє ь√ юяЁхфхы хь ёЄрэфрЁЄэ√ь юсЁрчюь). ╬сючэрўшь ўхЁхч $s(x,\la)$ Ёх°хэшх фшЇЇхЁхэЎшры№эюую єЁртэхэш 
$$
-(s^{[1]}(x))'-\sigma(x) s^{[1]}(x)-\sigma^2(x)s(x)=\la s(x)
$$
ё эрўры№э√ьш єёыютш ьш $s(0,\la)=0$, $s^{[1]}(0,\la)=\sqrt{\la}$ (Єръюх Ёх°хэшх ёє∙хёЄтєхЄ ш хфшэёЄтхээю, ёь.
\cite{SS2}). ╧єёЄ№, ъръ ш т яхЁтюь ярЁруЁрЇх, $\{\la_k\}_1^\infty$ --- ёюсёЄтхээ√х чэрўхэш  юяхЁрЄюЁр $L_D$ (эєыш
ЇєэъЎшш $s(\pi,\la)$), чрэєьхЁютрээ√х т яюЁ фъх тючЁрёЄрэш  ьюфєы . ┬ ёыєўрх тх∙хёЄтхээющ ЇєэъЎшш $\sigma$ юяхЁрЄюЁ
$L_D$ ёрьюёюяЁ цхэ ш тёх ёюсёЄтхээ√х чэрўхэш  $\la_k$ тх∙хёЄтхээ√ ш яЁюёЄ√. ╫хЁхч $\{\al_k\}_1^\infty$ ь√, яю-яЁхцэхьє,
юсючэрўрхь эюЁьшЁютюўэ√х ўшёыр, юяЁхфхыхээ√х т \eqref{eq:nc}.

═р°р Ўхы№ --- фрЄ№ юЎхэъш яюуЁх°эюёЄш ьхцфє ЇєэъЎшхщ $\sigma_N$, тюёёЄрэютыхээющ яю ъюэхўэюьє эрсюЁє ёяхъЄЁры№э√ї
фрээ√ї, ш шёЄшээющ ЇєэъЎшхщ $\sigma$. ─ы  ЇюЁьєышЁютъш Єръшї Ёхчєы№ЄрЄют эрь эхюсїюфшью ттхёЄш яЁюёЄЁрэёЄтр яюЄхэЎшрыют
(фрыхх фы  ъЁрЄъюёЄш ь√ эрч√трхь яюЄхэЎшрыюь ЇєэъЎш■ $\sigma$) ш яЁюёЄЁрэёЄтр ёяхъЄЁры№э√ї фрээ√ї. ╧юыюцшь
\begin{equation}\label{eq:sk}
    s_{2k-1} =\al_k-\frac{\pi}2, \quad s_{2k} =\sqrt{\lambda_k} - k,  \quad k= 1,2, \dots,  .
\end{equation}
▀ёэю, ўЄю яюёыхфютрЄхы№эюёЄ№ $\{s_k\}_1^\infty$ юсэючэрўэю юяЁхфхы хЄё  ш, юсЁрЄэю, юфэючэрўэю юяЁхфхы хЄ ёяхъЄЁры№э√х
фрээ√х $\{\al_k\}_1^\infty\cup\{\la_k\}_1^\infty$. ─рыхх ь√ сєфхь ЁрсюЄрЄ№ шьхээю ё яюёыхфютрЄхы№эюёЄ№■
$\{s_k\}_1^\infty$, ъюЄюЁє■ эрчютхь {\it Ёхуєы Ёшчютрээ√ьш ёяхъЄЁры№э√ьш фрээ√ьш}.

╥хяхЁ№  яюёЄЁюшь яЁюёЄЁрэёЄтр, ъюЄюЁ√ь  яЁшэрфыхцрЄ Ёхуєы Ёшчютрээ√х ёяхъЄЁры№э√х фрээ√х. ╬сючэрўшь ўхЁхч
$l^{\,\theta}_2$ тхёютюх $l_2$-яЁюёЄЁрэёЄтю, ёюёЄю ∙хх шч яюёыхфютрЄхы№эюёЄхщ ъюьяыхъёэ√ї ўшёхы $\bold x=\{x_1,
x_2,\dots\}$, Єръшї, ўЄю
$$
\|\bold x\|^2_\theta : =\sum_1^\infty |x_k|^2\, k^{2\theta} <\infty.
$$
╟рьхЄшь, ўЄю т ёыєўрх $\theta=1$ шч Ёрчыюцхэшщ \eqref{eq:spas} ш \eqref{eq:ncas} ёыхфєхЄ $\{\s_k\}_1^\infty\in
l^{\,\theta}_2\Longleftrightarrow q_0=0$, Є.х. яЁюёЄЁрэёЄтю ёяхъЄЁры№э√ї фрээ√ї эх ёютярфрхЄ ё $l_2^{1}$. ╤ єтхышўхэшхь
ярЁрьхЄЁр $\theta$ ёшЄєрЎш  єёыюцэ хЄё .
\begin{THEOREM}
╧єёЄ№ $\sigma(x)\in W_2^m$, Є.х. $q(x)\in W_2^{m-1}$, $m\geq1$. ╥юуфр фы  ёюсёЄтхээ√ї чэрўхэшщ $\la_k$ ш эюЁьшЁютюўэ√ї
ўшёхы $\al_k^2$ юяхЁрЄюЁр $L_D$ ёяЁртхфышт√ ёыхфє■∙шх ЇюЁьєы√. ╧Ёш эхўхЄэюь $m=2s+1$
\begin{align}\label{eq:laal1}
\sqrt{\la_k}&=k+\dfrac{h_0}{(2k)}+\dfrac{h_1}{(2k)^3}+\dots+\dfrac{h_s}{(2k)^{2s+1}}-
(-1)^s\dfrac{a_{2k}}{2(2k)^{2s+1}}+\dfrac{\gamma_{2k}}{k^{2s+2}},\\
\al_k&=\dfrac\pi2+\dfrac{g_1}{(2k)^2}+\dfrac{g_2}{(2k)^4}+\dots+ \dfrac{g_s}{(2k)^{2s}}- (-1)^s\dfrac{\wt
b_{2k}}{(2k)^{2s+1}}+ \dfrac{g_{s+1}}{(2k)^{2s+2}}+\dfrac{\gamma'_{2k}}{k^{2s+2}}.\notag
\end{align}
╧Ёш ўхЄэюь $m=2s$
\begin{align}\label{eq:laal2}
\sqrt{\la_k} &=k+\dfrac{h_0}{(2k)}+\dfrac{h_1}{(2k)^3}+\dots+\dfrac{h_{s-1}}{(2k)^{2s-1}}-
(-1)^s\dfrac{b_{2k}}{2(2k)^{2s}}+\dfrac{h_{s}}{(2k)^{2s+1}}+\dfrac{\gamma_{2k}}{k^{2s+1}},\\
\notag \al_k &=\dfrac\pi2+\dfrac{g_1}{(2k)^2}+\dfrac{g_2}{(2k)^4}+\dots+\dfrac{g_s}{(2k)^{2s}}- (-1)^s\dfrac{\wt
a_{2k-1}}{2(2k)^{2s}}+\dfrac{\gamma'_{2k-1}}{k^{2s+1}}.
\end{align}
┬ ¤Єшї ЇюЁьєырї $ \{\gamma_k\}$ ш $\{\gamma'_k\}$ --- яюёыхфютрЄхы№эюёЄш шч $l_2$ ш шї $l_2$-эюЁь√ юЎхэштр■Єё 
яюёЄю ээющ $C$,  чртшё ∙хщ юЄ $R$, эю эх чртшё ∙хщ юЄ $\sigma$  т °рЁх $\|\sigma\|_{m} \leq R$. ╫шёыр $a_l$, $\wt a_l$,
$b_l$ ш $\wt b_l$ (ь√ эх єърч√трхь шї чртшёшьюёЄ№ юЄ $m$)  юяЁхфхы ■Єё  ЇюЁьєырьш
$$
a_l=\frac2\pi\intl_0^\pi \sigma^{(m)}(t)\cos ltdt,\quad \wt a_l=\intl_0^\pi \left(\sigma (t)(\pi-t)\right)^{(m)}\,\cos
lt\, dt,
$$
$$
b_l=\frac2\pi\intl_0^\pi \sigma^{(m)}(t)\sin ltdt,\quad \wt b_l=\intl_0^\pi\left(\sigma (t)(\pi-t)\right)^{(m)} \sin
ltdt, \qquad l=1,2,\dots,
$$
р ўшёыр $h_j, g_j, \ 0\leq j\leq s+1$    ты ■Єё  эхяЁхЁ√тэ√ьш ш ЁртэюьхЁэю юуЁрэшўхээ√ьш т ърцфюь °рЁх ЇєэъЎшюэрырьш юЄ
$\sigma\in W_2^m$ ш шї ышэхщэ√х ўрёЄш $h^0_j,g^0_j$ яЁш $j\leq s-1$ т√Ёрцр■Єё  ЇюЁьєырьш
$$
h_j^0=(-1)^{j}\pi^{-1}\left(\sigma^{(2j)}(\pi)-\sigma^{(2j)}(0)\right),\quad g_j^0=(-1)^{j+1} \pi^{-1} \left(\sigma
(t)(\pi -t)\right)^{(2j-1)}\, \vert_0^\pi.
$$
┬ ёыєўрх $m=2s+1$ ¤Єш ЇюЁьєы√ ёюїЁрэ ■Єё  ш фы  $j=s$, р $g^0_{s+1}=0$. ┬ ёыєўрх  $m=2s$
$$
h^0_s=0,\qquad g_s^0=(-1)^{s+1} \pi^{-1} \left(\sigma (t)(\pi -t)\right)^{(2s-1)}\, \vert_0^\pi.
$$
┬ ўрёЄэюёЄш,
\begin{gather*}
h_0=\frac1{\pi}(\sigma(\pi)-\sigma(0)),\qquad g_1=h_0+\sigma'(0)+\pi^3h_0^2-\frac{\pi}2(\sigma^2(\pi)-\sigma^2(0)),\\
h_1=-\frac1{\pi}(\sigma''(\pi)-\sigma''(0))+\frac1{\pi}\int_0^\pi(\sigma'(x))^2\,dx-2h_0^2.
\end{gather*}
\end{THEOREM}
\begin{proof}   ╘юЁьєы√ \eqref{eq:laal1} ш \eqref{eq:laal2} ё юЎхэърьш $\{\gamma_k\}\in l_2$ ш $\{\gamma_k'\}\in l_2$ яЁшэрфыхцшЄ
┬.└.╠рЁўхэъю, ёь. \cite{Ma}. ┴юыхх ъюьяръЄэюх фюърчрЄхы№ёЄтю ё ЁртэюьхЁэ√ьш яю °рЁє $\|\sigma\|_m\le R$ юЎхэърьш
яюыєўхэю ртЄюЁрьш т \cite{SS4}.
\end{proof}
─ы  яюёЄЁюхэш  яЁюёЄЁрэёЄтр ёяхъЄЁры№э√ї фрээ√ї ь√ ЁрёёьюЄЁшь ёяхЎшры№э√х яюёыхфютрЄхы№эюёЄш
\begin{gather}\label{eq:e} \bold e_{2s-1} =
\{\, 0,\  2^{-(2s-1)},\, 0,\ 4^{-(2s-1)},\, 0,\ 6^{-(2s-1)},\ldots\}\qquad\text{ш}\\ \quad \bold e_{2s}= \{2^{-(2s)},\,
0,\  4^{-(2s)},\, 0,\  6^{-(2s)},\ldots\},
 \qquad s=1,2,\dots .
\end{gather}
╟рьхЄшь, ўЄю яюёыхфютрЄхы№эюёЄ№ $\bold e_p $ яЁшэрфыхцшЄ яЁюёЄЁрэёЄтє $l_2^{\,\theta}$ т ЄюўэюёЄш Єюуфр, ъюуфр
$0\leqslant\theta < p-1/2$. ─ы  ЇшъёшЁютрээюую $\theta\geqslant 0$ яюыюцшь $m=[\theta+1/2]$ (чфхё№ $[\cdot]$ --- Ўхыр 
ўрёЄ№ ўшёыр). ─ы  Єръюую $\theta$ юяЁхфхышь яЁюёЄЁрэёЄтю $\ell_D^{\,\theta}$ ъръ ъюэхўэюьхЁэюх Ёрё°шЁхэшх яЁюёЄЁрэёЄтр
$l_2^{\,\theta}$ ёыхфє■∙шь юсЁрчюь
$$
\ell_D^{\,\theta}=l_2^{\,\theta}\oplus \text{span}\{\bold e_k\}_{k=1}^{m}.
$$
╥ръшь юсЁрчюь,  $\ell_D^{\,\theta}$ ёюёЄюшЄ шч ¤ыхьхэЄют $\bold{x}+\sum_{k=1}^{m} c_k \bold e_k$, уфх $\bold x\in
l_2^{\,\theta}$, р $\{c_k\}_1^{m}$ --- яЁюшчтюы№э√х ъюьяыхъёэ√х ўшёыр. ╤ъры Ёэюх яЁюшчтхфхэшх ¤ыхьхэЄют шч
$l_D^{\,\theta}$ юяЁхфхы хЄё  ЇюЁьєыющ
$$
(\bold{x} +\sum_{k=1}^{m} c_k \bold e_k , \ \bold{y} +\sum_{k=1}^{m} d_k \bold e_k ) = (\bold x , \bold y)_\theta +
\sum_{k=1}^{m} c_k\overline{d_k},
$$
уфх $(\bold x , \bold y)_\theta$ --- ёъры Ёэюх яЁюшчтхфхэшх т $l_2^{\,\theta}$. ╧юёЄЁюхээюх яЁюёЄЁрэёЄтю  ёт цхь ё
Ёхуєы Ёшчютрээ√ьш ёяхъЄЁры№э√ьш фрээ√ьш фы  юяхЁрЄюЁр $L_D$. ╒юЄ  ¤Єю яЁюёЄЁрэёЄтю юяЁхфхыхэю ъръ ъюэхўэюьхЁэюх
Ёрё°шЁхэшх тхёютюую яЁюёЄЁрэёЄтр $l_2^{\,\theta}$, хую ¤ыхьхэЄ√ єфюсэхх чряшё√трЄ№ т ЇюЁьх юс√ўэ√ї яюёыхфютрЄхы№эюёЄхщ.
╧Ёш $\theta\in[0,1/2)$ яЁюёЄЁрэёЄтю  $\ell_D^{\,\theta}$ ёютярфрхЄ ё юс√ўэ√ь тхёют√ь яЁюёЄЁрэёЄтюь $l_2^{\,\theta}$.
╧Ёш $\theta\in[1/2,3/2)$ ¤ыхьхэЄ√ яЁюёЄЁрэёЄтр $\ell_D^{\,\theta}$ ёюёЄю Є шч яюёыхфютрЄхы№эюёЄхщ $\{x_k\}_1^\infty$ ё
ъююЁфшэрЄрьш
$$
x_k =y_k +\begin{cases} 0,\quad&
\text{хёыш }k\text{ эхўхЄэю,}\\
\alpha_1 \, k^{-1},\quad& \text{хёыш }k\text{ ўхЄэю,}
\end{cases}\quad \text{уфх}\ \ \{y_k\}_1^\infty \in
l_2^{\,\theta}, \ \ \alpha_1 \in \mathbb C,\ \text{ш Є.ф.}
$$
╦хуъю тшфхЄ№, ўЄю яЁюёЄЁрэёЄтю $\ell_D^{\,\eta}$  ъюьяръЄэю тыюцхэю т яЁюёЄЁрэёЄтю $\ell_D^{\,\theta}$ яЁш $\eta >
\theta$ (¤Єю ёЁрчє ёыхфєхЄ шч ъюьяръЄэюёЄш тыюцхэш  $l^{\,\eta}_2\hookrightarrow l_2^{\,\theta}$\ яЁш $\eta > \theta$).

╥хяхЁ№ юяЁхфхышь юЄюсЁрцхэшх $F$, шуЁр■∙хх ЎхэЄЁры№эє■ Ёюы№ т эр°хщ чрфрўх:
\begin{equation}\label{eq:F}
F(\sigma)=\{s_k\}_1^\infty,\qquad\sigma\in W_2^\theta/\{1\}.
\end{equation}
\begin{Theorem}\label{tm:2.1}
╧єёЄ№ $\{\la_k\}_1^\infty$ ш $\{\al_k\}_1^\infty$ -- яюёыхфютрЄхы№эюёЄш ёюсёЄтхээ√ї чэрўхэшщ ш эюЁьшЁютюўэ√ї ўшёхы
юяхЁрЄюЁр $L_D$ ё ъюьяыхъёэ√ь яюЄхэЎшрыюь $q(x)\in W_2^{\theta-1}$ (шыш $\sigma(x)\in W_2^\theta$), $\theta
>0$. ╧єёЄ№ яюёыхфютрЄхы№эюёЄ№ $\{s_k\}_{k\in\bZ_0}$, юяЁхфхыхэр
ЁртхэёЄтрьш \eqref{eq:sk}. ╥юуфр $F$  ты хЄё  юЄюсЁрцхэшхь шч $W_2^\theta$ т $\ell_D^{\,\theta}$, яЁшўхь
$$
F(\sigma)=S\sigma-\Phi(\sigma),
$$
уфх ышэхщэ√щ юяхЁрЄюЁ $S$ юяЁхфхыхэ ЁртхэёЄтрьш
\begin{equation*}
(S\sigma)_{2k}=-\frac1{\pi}\intl_0^\pi\sigma(t)\sin(2kt)dt,\qquad
(S\sigma)_{2k-1}=-\intl_0^\pi(\pi-t)\sigma(t)\cos(2kt)\, dt ,\quad k\in\bN,
\end{equation*}
ш  ты хЄё  єэшЄрЁэ√ь шчюьюЁЇшчьюь шч $W_2^\theta$ т $\ell_D^\theta$, р эхышэхщэюх юЄюсЁрцхэшх $\Phi$ юЄюсЁрцрхЄ
$W_2^\theta$ т $\ell_2^{\,\tau}$, уфх
$$
\tau=\left\{ \begin{array}{ll} 2\theta,\ \text{хёыш}\ 0<\theta\leqslant1,\\ \theta+1,\ \text{хёыш}\
1\leqslant\theta<\infty.\end{array} \right.
$$
╧Ёш ¤Єюь юЄюсЁрцхэшх $\Phi:\,W_2^\theta\to\ell_2^\tau$  ты хЄё  юуЁрэшўхээ√ь т ърцфюь °рЁх, Є.х.
\begin{equation*}
\|\Phi(\sigma)\|_\tau\leq C\|\sigma\|_\theta,
\end{equation*}
уфх яюёЄю ээр  $C$ чртшёшЄ юЄ $R$, эю эх чртшёшЄ юЄ $\sigma$ т °рЁх $\|\sigma\|_\theta\leq R$.
\end{Theorem}
─юърчрЄхы№ёЄтю ¤Єющ ЄхюЁхь√ ёь. т ЁрсюЄх ртЄюЁют \cite{SS4}.

─ы  ЇюЁьєышЁютъш Ёхчєы№ЄрЄют ттхфхь эхъюЄюЁ√х юсючэрўхэш . ╬сючэрўшь ўхЁхч  $W^{\,\theta}_{2,\mathbb R}$  ьэюцхёЄтю
тёхї тх∙хёЄтхээ√ї ЇєэъЎшщ шч $W^{\,\theta}_2$. ╫хЁхч $\Gamma^{\,\theta}$ юсючэрўшь  ьэюцхёЄтю  ЇєэъЎшщ $\sigma\in
W^{\,\theta}_{2,\mathbb R} /\{1\}$, фы  ъюЄюЁ√ї $\la_k(\sigma)\geqslant  k^2$, р ўхЁхч $\mathcal B^{\,\theta}_\Gamma
(R)$ --- яхЁхёхўхэшх ьэюцхёЄтр $\Gamma^{\,\theta}$ ё чрьъэєЄ√ь °рЁюь $\mathcal B^{\,\theta}_{\mathbb R}(R)$ Ёрфшєёр $R$
т яЁюёЄЁрэёЄтх $W^{\,\theta}_{2, \mathbb R}$.

┼ёыш $\sigma\in\Gamma^{\,\theta}$, Єю ёюсёЄтхээ√х чэрўхэш  юяхЁрЄюЁр $L_D$ яюфўшэхэ√ єёыютш ь $\lambda_1 <\la_2<\dots$,
$\la_k\ge k^2$. ─ы  Ёхуєы Ёшчютрээ√ї ёяхъЄЁры№э√ї фрээ√ї ¤Єш эхЁртхэёЄтр ¤ътштрыхэЄэ√ ёыхфє■∙шь
\begin{equation}\label{s1}
s_{2k} \geqslant 0, \qquad s_{2k}-s_{2k+2}<1, \qquad k=1,2,\dots .
\end{equation}
╙ёыютш  эхюЄЁшЎрЄхы№эюёЄш тёхї эюЁьшЁютюўэ√ї ўшёхы ¤ътштрыхэЄэ√ єёыютш ь
\begin{equation}\label{s2}
s_{2k-1}>-\pi/2,\qquad k=1,2,\dots.
\end{equation}
╧юёыхфютрЄхы№эюёЄ№ $\{s_k\}_1^\infty \in l_2$, яю¤Єюьє  фы  ы■сющ тх∙хёЄтхээющ ЇєэъЎшш $\sigma\in \Gamma^{\,\theta}$
эрщфхЄё  ўшёыю $h=h(\sigma)>0$, Єръюх, ўЄю
\begin{equation}\label{sh}
s_{2k} \geqslant 0, \qquad s_{2k}-s_{2k+2}\leqslant 1-h,\qquad s_{2k-1}\geqslant -\pi/2+h, \qquad k=1,2,\dots .
\end{equation}
╘шъёшЁєхь яЁюшчтюы№э√х ўшёыр $r>0$ ш $h \in (0,1)$. ╬сючэрўшь ўхЁхч $\Omega^{\,\theta}(r,h)$ ёютюъєяэюёЄ№ тх∙хёЄтхээ√ї
яюёыхфютрЄхы№эюёЄхщ $\{s_k\}_1^\infty$, фы  ъюЄюЁ√ї т√яюыэхэ√ эхЁртхэёЄтр \eqref{sh} ш ъюЄюЁ√х ыхцрЄ т чрьъэєЄюь °рЁх
Ёрфшєёр $r$ яЁюёЄЁрэёЄтр $\ell^{\,\theta}_D$,  Є.х. $\|\{s_k\}\|_\theta \leqslant r$  (чфхё№ ш фрыхх яюфЁрчєьхтрхь, ўЄю
$\|\cdot\|_\theta$ ючэрўрхЄ эюЁьє т яЁюёЄЁрэёЄтх  $\ell^{\,\theta}_D$). ╫хЁхч
$\Omega^{\,\theta}=\Omega^{\,\theta}(\infty,0)$ юсючэрўшь ьэюцхёЄтю тёхї тх∙хёЄтхээ√ї яюёыхфютрЄхы№эюёЄхщ
$\{s_k\}_1^\infty \in l_D^{\,\theta}$, фы  ъюЄюЁ√ї ёяЁртхфышт√ эхЁртхэёЄтр \eqref{s1} ш \eqref{s2}. ╬сючэрўшь ўхЁхч
$\widehat \Omega^{\,\theta}$ ьэюцхёЄтю яюёыхфютрЄхы№эюёЄхщ $\{s_k\}_1^\infty \in \ell^{\,\theta}_D$, фы  ъюЄюЁ√ї ўшёыр
$\lambda_k =(s_{2k} +k)^2$ юсЁрчє■Є ёЄЁюую тючЁрёЄр■∙є■ яюёыхфютрЄхы№эюёЄ№ (яЁш ¤Єюь т ёыєўрх $\lambda_k <0$ ъююЁфшэрЄ√
$\{s_{2k}\}$  сєфєЄ эхтх∙хёЄтхээ√ьш), р тёх ўшёыр $\alpha_k = s_{2k-1}+\pi/2$ яюыюцшЄхы№э√.
\begin{Theorem}\label{tm:2.2} ╧єёЄ№ $\theta\ge0$ ЇшъёшЁютрэю, р $F$ --- юЄюсЁрцхэшх, ттхфхээюх т \eqref{eq:F}.\hfill\break
{\bf 1)}\quad ╬ЄюсЁрцхэшх \ $F: \Gamma^{\,\theta} \to \Omega^{\,\theta}$\ хёЄ№ сшхъЎш .

\noindent{\bf 2)}\quad ╬ЄюсЁрцхэшх $F: W^\theta_{2,\mathbb R}/\{1\} \to \widehat\Omega^{\,\theta}$ Єръцх хёЄ№ сшхъЎш .
╚э√ьш ёыютрьш, ўшёыр $\{\lambda_k\}^\infty_1$ ш $\{\alpha_k\}^\infty_1$ яЁхфёЄрты ■Є ёяхъЄЁ ш эюЁьшЁютюўэ√х ўшёыр
юяхЁрЄюЁр $L_D$  хёыш ш Єюы№ъю хёыш яхЁтр  яюёыхфютрЄхы№эюёЄ№  ты хЄё  ёЄЁюую ьюэюЄюээющ,  тЄюЁр  ёюёЄюшЄ шч
яюыюцшЄхы№э√ї ўшёхы, р юсЁрчютрээр  шч эшї яю ЇюЁьєырь \eqref{eq:sk} яюёыхфютрЄхы№эюёЄ№ $\{s_k\}^\infty_1$ яЁшэрфыхцшЄ
яЁюёЄЁрэёЄтє $\ell^{\,\theta}_D$.

\noindent{\bf 3)}\quad ╧єёЄ№ $R$ яЁюшчтюы№эюх яюыюцшЄхы№эюх ўшёыю. ═рщфєЄё  яюыюцшЄхы№э√х ўшёыр $r=r(R), h=h(R)$,
Єръшх, ўЄю
$$
F (\mathcal B^{\,\theta}_\Gamma (R) ) \subset \Omega^\theta (r,h).
$$
\noindent{\bf 4)}\quad ╤яЁртхфыштю юсЁрЄэюх єЄтхЁцфхэшх: фы  ы■с√ї ўшёхы $r>0$ ш $h\in(0,1)$ эрщфхЄё  ўшёыю $R>0$,
Єръюх, ўЄю
$$
F^{-1}(\,\Omega^{\,\theta}(r,h)) \subset \mathcal B^{\,\theta}_\Gamma (R).
$$
\noindent{\bf 5)}\quad ╤яЁртхфыштю яЁхфёЄртыхэшх
$$
F^{-1} = S^{-1} +\Psi, \qquad \Psi: \Omega^{\,\theta} \to W^\tau_2,
$$
уфх ўшёыю $\tau$  юяЁхфхыхэю т ╥хюЁхьх \ref{tm:2.1}. ╬ЄюсЁрцхэшх $\Phi: \Omega^\theta \to W^\tau_2,$ рэрышЄшўэю, яЁшўхь
\begin{equation*}
\|\Psi\bold y\|_\tau \leqslant C\|\bold y\|_\theta \qquad\text{фы  тёхї} \ \, \bold y \in\Omega^{\,\theta}(r,h),
\end{equation*}
уфх яюёЄю ээр  $C$  чртшёшЄ Єюы№ъю юЄ $r$  ш $h$.

\noindent{\bf 6)}\quad ╧єёЄ№ яюёыхфютрЄхы№эюёЄш $\bold y, \bold y_1$ Ёхуєы Ёшчютрээ√ї ёяхъЄЁры№э√ї фрээ√ї ыхцрЄ т
ьэюцхёЄтх $\Omega^{\,\theta}(r,h)$. ╥юуфр  яЁююсЁрч√ $\sigma = F^{-1}\bold y,\ \,\sigma_1 = F^{-1}\bold y_1$  ыхцрЄ т
ьэюцхёЄтх $\mathcal B_\Gamma^{\,\theta} (R)$ ш ёяЁртхфышт√ юЎхэъш
\begin{equation}\label{a}
C_1\|\bold y -\bold y_1\|_\theta \leqslant \|\sigma -\sigma_1\|_\theta \leqslant C_2\|\bold y -\bold y_1\|_\theta,
\end{equation}
уфх ўшёыю $R$ ш яюёЄю ээ√х $C_1, C_2$ чртшё Є Єюы№ъю $r$ ш $h$.

\noindent{\bf 7)}\quad ╬сЁрЄэю, хёыш $\sigma, \sigma_1$ ыхцрЄ т °рЁх $\mathcal B_\mathbb R^{\,\theta}(R)$, Єю
яюёыхфютрЄхы№эюёЄш $\bold y, \bold y_1$ Ёхуєы Ёшчютрээ√ї ёяхъЄЁры№э√ї фрээ√ї ¤Єшї ЇєэъЎшщ ыхцрЄ т ьэюцхёЄтх
$\Omega^\theta(r,h)$  ш ёяЁртхфышт√ юЎхэъш
\begin{equation}\label{b}
C_1 \|\sigma -\sigma_1\|_\theta \leqslant\|\bold y -\bold y_1\|_\theta \leqslant C_2 \|\sigma -\sigma_1\|_\theta .
\end{equation}
╟фхё№ ўшёыр $r>0, \ \, h\in (0, 1)$  ш яюёЄю ээ√х $C_1$ ш $C_2$ чртшё Є Єюы№ъю юЄ $R$.
\end{Theorem}
\begin{proof}
╤ыєўрщ $\theta>0$ ЁрёёьюЄЁхэ т ёЄрЄ№х ртЄюЁют \cite{SS6}. ╤ыєўрщ $\theta=0$ шчєўхэ т ЁрсюЄх ╨.~╬.~├Ёшэштр \cite{Hr}.
\end{proof}
╫шёыю $R$   ш яюёЄю ээ√х $C_2, C^{-1}_1$ т \eqref{a} єтхышўштр■Єё  яЁш $r\to\infty$  шыш $h\to 0$. ╫шёыр $ r,\, h^{-1},
C_2 $ ш $C^{-1}_1$ т \eqref{b} Єръцх єтхышўштр■Єё  яЁш $R\to\infty$. ╬Ўхэъш эр ЁюёЄ ¤Єшї ъюэёЄрэЄ ртЄюЁрьш яюыєўхэ√, эю
фюърчрЄхы№ёЄтр ёююЄтхЄёЄтє■∙шї єЄтхЁцфхэшщ ЄЁхсє■Є чэрўшЄхы№эющ ЁрсюЄ√, ш ртЄюЁ√ эрьхЁхтр■Єё  яюёт ЄшЄ№ ¤Єюьє юЄфхы№эє■
ёЄрЄ№■.

╤ яюью∙№■ Ёхчєы№ЄрЄют ЄхюЁхь√ \ref{tm:2.2} ь√ яюыєўшь ЄхяхЁ№ юЎхэъш яюуЁх°эюёЄш т чрфрўх тюёёЄрэютыхэш  яюЄхэЎшрыр яю
ъюэхўэюьє эрсюЁє ёяхъЄЁры№э√ї фрээ√ї.  ╩ръ ш т яхЁтюь ярЁруЁрЇх ь√ яЁхфырурхь фтх яюёЄрэютъш ¤Єющ чрфрўш: яЁш эрышўшш
ряЁшюЁэющ шэЇюЁьрЎшш ю яюЄхэЎшрых ш яЁш эрышўшш ряЁшюЁэющ шэЇюЁьрЎшш ю ёяхъЄЁры№э√ї фрээ√ї. ╠√ эрўэхь ё яхЁтюую ёыєўр .
╧Ёхфяюыюцшь, ўЄю эрь шчтхёЄэр ряЁшюЁэр  шэЇюЁьрЎш  $\sigma\in W_{2,\bR}^\theta/\{1\}$ фы  эхъюЄюЁюую ЇшъёшЁютрээюую
$\theta>0$ ш $\|\sigma\|_\theta\le R$ фы  эхъюЄюЁюую $R>0$ (т эр°шї юсючэрўхэш ї $\sigma\in B^\theta_{\bR}(R)$). ╧єёЄ№
$\{\la_k\}_1^\infty$
--- ёюсёЄтхээ√х чэрўхэш , р $\{\al_k\}_1^\infty$ --- эюЁьшЁютюўэ√х ўшёыр юяхЁрЄюЁр $L_D$ ё яюЄхэЎшрыюь $\sigma$.
╬сючэрўшь ўхЁхч $\bold
s=\{s_k\}_1^\infty=F(\sigma)$ Ёхуєы Ёшчютрээ√х ёяхъЄЁры№э√х фрээ√х
--- ёюуырёэю ЄхюЁхьх \ref{tm:2.1}, $\bold s\in\ell^\theta_D$. ▌Єю ючэрўрхЄ, ўЄю
\begin{equation}\label{eq:sd}
\bold s=c_1\bold e_1+c_2\bold e_2+\dots+c_m\bold e_m+\bold a,
\end{equation}
уфх $m=[\theta+1/2]$, яюёыхфютрЄхы№эюёЄш $\bold e_k$ юяЁхфхыхэ√ т \eqref{eq:e}, р $\bold a=\{a_k\}_1^\infty\in
l_2^\theta$. ╧Ёхфяюыюцшь, ўЄю эрьш ё яюуЁх°эюёЄ№■ $\eps$ т√ўшёыхэ√ $N$ ёюсёЄтхээ√ї чэрўхэшщ $\{\wt\la_k\}_1^N$ ш $N$
эюЁьшЁютюўэ√ї ўшёхы $\{\wt\al_k\}_1^N$, р шьхээю, т√яюыэхэ√ эхЁртхэёЄтр \eqref{eq:err1}. ╧Ёхфяюыюцшь Єръцх, ўЄю эрь
шчтхёЄэ√ ўшёыр $\wt c_j$ шч яЁхфёЄртыхэш  \eqref{eq:sd} фы  тёхї $1\le j\le m$. {\it ╩юэхўэ√ь эрсюЁюь яЁшсышцхээ√ї
ёяхъЄЁры№э√ї фрээ√ї} ь√ эрчютхь
$$
\{c_j\}_1^m\cup\{\wt\la_k\}_1^N\cup\{\wt\al_k\}_1^N.
$$
╥хяхЁ№ яю ёяхъЄЁры№э√ь фрээ√ь ь√ эрщфхь $2N$--ряяЁюъёшьрЎш■ яюЄхэЎшрыр. ╧юёЄЁюшь тэрўрых тёяюьюурЄхы№эє■ ЇєэъЎш■
$\sigma_0$
--- яЁюшчтюы№эє■ ЇєэъЎш■ шч $W_2^\theta$ ё єёыютшхь
$$
\left(F(\sigma_0)-\sum_{j=1}^mc_j\bold e_j\right)\in l_2^\theta.
$$
┬ ўрёЄэюёЄш, яЁш $\theta\in[0,1/2)$ ўшёыю $m=0$ ш $\sigma_0=0$. ╧Ёш $\theta\in[1/2,3/2)$ ўшёыю $m=1$ ш ЇєэъЎш■
$\sigma_0$ ьюцэю т√сЁрЄ№ Ёртэющ $c_1x$ --- ыхуъю тшфхЄ№, ўЄю фы  Єръющ ЇєэъЎшш $\al_k=\pi/2$, р $\la_k=k^2+c_1$. ╧Ёш
$\theta\in[3/2,5/2)$ ьюцэю, эряЁшьхЁ, тч Є№
$$
\sigma_0(x)=c_1x+(c_2-2c_1-\pi^3c_1^2/2)x(\pi-x).
$$
╚ч яЁштхфхээ√ї т ЄхюЁхьх \ref{tm:3} ЇюЁьєы фы  ЇєэъЎшюэрыют $h_0$ ш $g_1$ ёыхфєхЄ, ўЄю фы  Єръющ ЇєэъЎшш $h_0=c_1$, р
$g_1=c_2$. ╧Ёш $\theta\in[5/2,7/2)$
$$
\sigma_0(x)=\al x+\beta x(\pi-x)+\gamma x^2(\pi-x),
$$
уфх ъю¤ЇЇшЎшхэЄ√ $\al$, $\beta$, $\gamma$ т√Ёрцр■Єё  ўхЁхч $c_1$, $c_2$ ш $c_3$ т  тэюь тшфх. ╠√ эх яЁштюфшь чфхё№
уЁюьючфъшх ЇюЁьєы√, яюёъюы№ъє тшф ЇєэъЎшш $\sigma_0$ эх тыш хЄ эр фры№эхщ°шх Ёрёёєцфхэш . ╧Ёш сюы№°шї $\theta$ эрщЄш
 тэ√щ рэрышЄшўхёъшщ тшф ЇєэъЎшш $\sigma_0$ ёыюцэхх
--- фы  ¤Єюую ЄЁхсєхЄё   тэю т√яшёрЄ№ ЇєэъЎшюэры√ $h_j$ ш $g_j$ т Ёрчыюцхэш ї \eqref{eq:laal1} ш \eqref{eq:laal2}. ═рь трцэю
ыш°№ Єю, ўЄю ЇєэъЎш■ $\sigma_0$ тёхуфр ьюцэю т√сЁрЄ№ Єръ, ўЄю $\|\sigma_0\|_\theta\le C(R)$. ─ы  ¤Єюую фюёЄрЄюўэю
яюыюцшЄ№ $\sigma_0=F^{-1}(\bold e)$ (фы  ъЁрЄъюёЄш ь√ юсючэрўшыш $\bold e=\sum_{j=1}^mc_j\bold e_j$) ш тюёяюы№чютрЄ№ё 
эхЁртхэёЄтрьш \eqref{b}:
$$
\|\sigma_0\|_\theta\le C_1^{-1}\|\bold e\|_\theta\le C_1^{-1}\|\bold s\|_\theta\le C_2C_1^{-1}\|\sigma\|_\theta.
$$
╚Єръ, ЇєэъЎш  $\sigma_0$ яюёЄЁюхэр. ╫хЁхч $\la_k^0$ ш $\al_k^0$ ь√ юсючэрўшь ёяхъЄЁры№э√х фрээ√х, р ўхЁхч $\bold
s^0=\bold e+\bold a^0$ юсючэрўшь Ёхуєы Ёшчютрээ√х ёяхъЄЁры№э√х фрээ√х фы  яюЄхэЎшрыр $\sigma_0$. ╥хяхЁ№, шёяюы№чє 
рыуюЁшЄь, юяшёрээ√щ т ЁрсюЄх \cite{SS6} (юэ тюёїюфшЄ ъ ЁрсюЄрь ╧ю°хы  ш ╥ЁєсютшЎр, ёь. \cite{PT}) яюёЄЁюшь яюЄхэЎшры
$\wt\sigma_N$, ёяхъЄЁры№э√х фрээ√х ъюЄюЁюую шьх■Є тшф
\begin{equation}\label{eq:approx}
\{c_j\}_1^m\cup\{\wt\la_1,\dots,\wt\la_N,\la_{N+1}^0,\la_{N+2}^0\}_1^N\cup\{\wt\al_1,\dots,\wt\al_N,\al_N^0,\al_{N+1}^0,\dots\}_1^N.
\end{equation}
╘єэъЎш■ $\wt\sigma_N$ ь√ эрчютхь {\it $2N$--ряяЁюъёшьрЎшхщ яюЄхэЎшрыр $\sigma$.} ┼ёЄхёЄтхээю, фрээ√х шчьхЁхэшщ фюыцэ√
с√Є№ ёюуырёютрэ√ ё ряЁшюЁэющ шэЇюЁьрЎшхщ. ┬ ърўхёЄтх Єръюую єёыютш  ёюуырёютрэш  ь√ яюЄЁхсєхь
$\|\wt\sigma_N\|_\theta\le R$.
\begin{Theorem}\label{tm:main}
╧єёЄ№  $\sigma\in B^\theta_{\Gamma}(R)$ фы  эхъюЄюЁюую $\theta>0$, р $\wt\sigma_N$ --- хх $2N$--ряяЁюъёшьрЎш ,
юяЁхфхыхээр  т√°х, яЁшўхь $\|\wt\sigma_N\|_\theta\le R$. ╥юуфр фы  ы■сюую $\tau\in[0,\theta)$ т√яюыэхэ√ юЎхэъш
\begin{equation}\label{eq:main}
\|\sigma-\wt\sigma_N\|_\tau\le C_1^{-1}\|\bold s-\bold{\wt s}\|_\tau\le C(R)N^{\tau-\theta}+\eps
C(R)\begin{cases}N^{\tau -1/2}\ \ &\text{хёыш} \ \, \tau >1/2,\\
(\ln N)^{1/2} \ \ &\text{хёыш} \ \,  \tau = 1/2,\\
 1, \ \ &\text{хёыш}  \ \, \tau < 1/2.\end{cases}
\end{equation}
┬ ўрёЄэюёЄш,
\begin{equation}\label{eq:main1}
\|\sigma(x) - \wt\sigma_{N}(x)\|_\tau \leqslant C(R)\left\{\begin{aligned} &\varepsilon^\gamma, \ \gamma
=\frac{2(\theta-\tau)}{2\theta-1},\ \ \ & &\text{хёыш} \ \,
\tau > 1/2,\ \, &&\text{р} \ \, N= \varepsilon^{-2/(2\theta-1)},\\
&\varepsilon|\ln\varepsilon|^{1/2}, \ \ \ & &\text{хёыш} \ \, \tau =1/2, \ \, &&\text{р} \ \, N=\varepsilon^{-1/(\theta -\tau)},\\
&\varepsilon, \ \ \ &&\text{хёыш} \ \, \tau < 1/2, \ \, &&\text{р} \ \, N\geqslant \varepsilon^{-1/(\theta -\tau)}.
\end{aligned}\right.
\end{equation}
╩юэёЄрэЄ√ т \eqref{eq:main} ш \eqref{eq:main1} чртшё Є Єюы№ъю юЄ $R$, $\theta$ ш $\tau$.
\end{Theorem}
\begin{proof}
╧єёЄ№
$$
\bold {\wt s}=\{\wt s_k=s_k(\sigma_N)\}_1^\infty, \qquad\bold{\wt s}=c_1\bold e_1+\dots+c_m\bold e_m+\bold{\wt a}
$$
--- Ёхуєы Ёшчютрээ√х ёяхъЄЁры№э√х фрээ√х, яюёЄЁюхээ√х яю ЇєэъЎшш $\sigma_N$.
╥юуфр
$$
\|\bold s-\bold{\wt s}\|_\tau^2=\sum_{k=1}^{\infty}|a_k-\wt a_k|^2k^{2\tau}\le
\eps^2\sum_{k=1}^{2N}k^{2\tau-2}+\sum_{k=2N+1}^\infty|a_k-a_k^0|^2k^{2\tau}.
$$
╧хЁтр  ёєььр яЁш $\tau>1/2$ юЎхэштрхЄё  тхышўшэющ $C\eps^2N^{2\tau-1}$ ё рсёюы■Єэющ ъюэёЄрэЄющ $C$. ╧Ёш $\tau=1/2$ ¤Єр
ёєььр юЎхэштрхЄё  тхышўшэющ $C\eps^2\ln N$, р яЁш $\tau\in[0,1/2)$ фюяєёърхЄ юЎхэъє тхышўшэющ $C\eps^2$. ┬ЄюЁє■ ёєььє
юЎхэшь, яЁшьхэшт эхЁртхэёЄтю \eqref{b}
\begin{multline*}
\sum_{k=2N+1}^\infty|a_k-a_k^0|^2k^{2\tau}\le2\sum_{k=2N+1}^\infty
k^{2\tau-2\theta}\left(|a_k|^2+|a_k^0|^2\right)k^{2\theta}\le\\
\le2(2N)^{2\tau-2\theta}\left(\|\bold a\|_\theta^2+\|\bold a^0\|_\theta^2\right)\le
2C_2\left(\|\sigma\|_\theta^2+\|\sigma_0\|_\theta^2\right)(2N)^{2\tau-2\theta}\le C(R)N^{2\tau-2\theta}.
\end{multline*}
╚Єръ,
\begin{equation*}
\|\bold s-\bold{\wt s}\|_\tau\le C(R)N^{\tau-\theta}+\eps
C(R)\begin{cases}N^{\tau -1/2}\ \ &\text{хёыш} \ \, \tau >1/2,\\
(\ln N)^{1/2} \ \ &\text{хёыш} \ \,  \tau = 1/2,\\
 1, \ \ &\text{хёыш}  \ \, \tau < 1/2,\end{cases}.
\end{equation*}
┬эют№ яЁшьхэ   \eqref{b}, яЁшфхь ъ юЎхэъх \eqref{eq:main}. ╬Ўхэър \eqref{eq:main1} эхьхфыхээю ёыхфєхЄ шч
\eqref{eq:main}.
\end{proof}
╨рёёьюЄЁшь ЄхяхЁ№ ёыєўрщ, ъюуфр т ърўхёЄтх ряЁшюЁэющ шэЇюЁьрЎшш т√сЁрэю єёыютшх эр ёяхъЄЁры№э√х фрээ√х. ╧Ёхфяюыюцшь,
ўЄю эрь шчтхёЄхэ ъюэхўэ√щ эрсюЁ ёяхъЄЁры№э√ї фрээ√ї $\{\wt\la_k\}_1^N\cup\{\wt\al_k\}_1^N$, шчьхЁхээ√ї ё ЄюўэюёЄ№■
$\eps$: яєёЄ№ т√яюыэхэ√ эхЁртхэёЄтр \eqref{eq:err1}. ╠√ сєфхь яЁхфяюырурЄ№, ўЄю фы  ¤Єюую эрсюЁр фрээ√ї т√яюыэхэ√
эхЁртхэёЄтр $\wt\la_k>k^2$, $k=1,\,\dots,\,N$. ▌Єш єёыютш  эх юуЁрэшўштр■Є юс∙эюёЄш, Є.ъ. ь√ тёхуфр ьюцхь ёфхырЄ№
яЁхюсЁрчютрэшх ёяхъЄЁры№э√ї фрээ√ї
$$
\wt\la_k\longrightarrow\wt\la_k+c,\qquad\al_k\longrightarrow\wt\al_k\cdot\sqrt{\frac{\wt\la_k+c}{\wt\la_k}},
$$
тюёёЄрэютшЄ№ яюЄхэЎшры $\sigma+c$ яю эютюьє эрсюЁє ёяхъЄЁры№э√ї фрээ√ї, р чрЄхь т√ўхёЄ№ ъюэёЄрэЄє $c$. ╩Ёюьх Єюую, ь√
сєфхь яЁхфяюырурЄ№ эрышўшх ряЁшюЁэющ шэЇюЁьрЎшш, р шьхээю, яЁхфяюыюцшь, ўЄю Ёхуєы Ёшчютрээ√х ёяхъЄЁры№э√х фрээ√х
шёъюьюую яюЄхэЎшрыр $\bold s=\{s_k\}_1^\infty$ ыхцрЄ т ьэюцхёЄтх $\Omega^\theta(r,h)$ фы  эхъюЄюЁюую $\theta>0$, $r>0$
ш $h\in(0,1)$. ╫шёыр $c_j$ шч яЁхфёЄртыхэш  \eqref{eq:sd} ь√, ъръ ш Ёрэхх, сєфхь ёўшЄрЄ№ шчтхёЄэ√ьш. ╘єэъЎш■ $\sigma_0$
ш $2N$--ряяЁюъёшьрЎш■ $\wt\sigma_N$ ь√ яюёЄЁюшь Єръ цх, ъръ ш Ёрэхх. ═ръюэхЎ, ь√ яЁхфяюыюцшь, ўЄю ряЁшюЁэр  шэЇюЁьрЎш 
ш Ёхчєы№ЄрЄ√ шчьхЁхэшщ ёюуырёютрэ√: фы  тхъЄюЁр $\bold{\wt s}$ Ёхуєы Ёшчютрээ√ї ёяхъЄЁры№э√ї фрээ√ї, яюёЄЁюхээ√ї яю
эрсюЁє \eqref{eq:approx}, т√яюыэхэю єёыютшх $\bold{\wt s}\in\Omega^\theta(r,h)$.
\begin{Theorem}
╧єёЄ№ юср тхъЄюЁр Ёхуєы Ёшчютрээ√ї ёяхъЄЁры№э√ї фрээ√ї $\bold s$ ш $\bold{\wt s}$ ыхцрЄ т $\Omega^\theta(r,h)$ фы 
эхъюЄюЁ√ї $\theta>0$, $r>0$ ш $h\in(0,1)$. ╧єёЄ№ $\wt\sigma_N$ --- $2N$--ряяЁюъёшьрЎш , юяЁхфхыхээр  т√°х. ╥юуфр фы 
ы■сюую $\tau\in[0,\theta)$ т√яюыэхэ√ юЎхэъш \eqref{eq:main} ш \eqref{eq:main1} ё чрьхэющ тёхї  ъюэёЄрэЄ $C(R)$ эр
тхышўшэ√ $C(r,h)$, чртшё ∙шх Єюы№ъю юЄ $\theta$, $\tau$, $r$ ш $h$.
\end{Theorem}
\begin{proof}
╩ръ ш т ЄхюЁхьх \ref{tm:main} юЎхэшь
$$
\|\bold s-\bold{\wt s}\|_\tau^2\le\sum_{k=2N+1}^\infty|a_k-a_k^0|^2k^{2\tau}+\eps
C\begin{cases}N^{\tau -1/2}\ \ &\text{хёыш} \ \, \tau >1/2,\\
(\ln N)^{1/2} \ \ &\text{хёыш} \ \,  \tau = 1/2,\\
 1, \ \ &\text{хёыш}  \ \, \tau < 1/2,\end{cases}
$$
ё рсёюы■Єэющ ъюэёЄрэЄющ $C$ (ўшёыр $\theta$ ш $\tau$ ь√ ёўшЄрхь ЇшъёшЁютрээ√ьш). ─рыхх,
\begin{multline*}
\sum_{k=2N+1}^\infty|a_k-a_k^0|^2k^{2\tau}\le2\sum_{k=2N+1}^\infty
k^{2\tau-2\theta}\left(|a_k|^2+|a_k^0|^2\right)k^{2\theta} \le\\ \le2\left(\|\bold a\|_\theta^2+\|\bold
a^0\|_\theta^2\right)(2N)^{2\tau-2\theta}\le C(r)N^{2\tau-2\theta}.
\end{multline*}
╧Ёшьхэ   ЄхяхЁ№ \eqref{a}, яЁшфхь ъ \eqref{eq:main} (ё чрьхэющ $C(R)$ эр $C(r,h)$). ╬Ўхэър \eqref{eq:main1} тэют№
ёыхфєхЄ шч \eqref{eq:main}.
\end{proof}

Savchuk A.M.

Moscow State University
e-mail: artem savchuk@mail.ru



\end{document}